\documentstyle[12pt]{article}%

\input amssym

\begin{document}
\vspace*{1cm}
\begin{center}

{\Large\bf Finite Groups with Some s-semipermutable subgroups\footnote{ Supported by the  major project of Basic and
Applied Research ( Natural Science)  in Guangdong Province£¬ China ( Grant Number:  2017KZDXM058) and
the Science and Technology Program of Guangzhou Municipality, China ( Grant number: 201804010088).}}\\[1cm]

\vskip0.5cm \centerline {Yangming LI\footnote{Corresponding
author.}}

\centerline {Dept. of Math., Guangdong University  of Education,
Guangzhou, 510310, China}
 \centerline {Email: liyangming@gdei.edu.cn}

\end{center}
\vspace{0.5cm}

 \begin{center}
\begin{minipage}{12cm}
 {\bf Abstract:}\hspace{0.3cm} Suppose that $G$ is a finite group and $H$ is a subgroup of $G$. We
say that $H$ is  s-semipermutable  in $G$ if $HG_p = G_pH$ for any
Sylow $p$-subgroup  $G_p$ of $G$ with $(p, |H|) = 1$.
 We investigate the influence of s-semipermutable subgroups on
the structure of finite groups. Some recent results are
generalized.\\

{\bf MSC(2000):}\hspace{0.5cm}20D10, 20D15

{\bf Keywords:}\hspace{0.5cm}  s-semipermutable subgroup, the smallest generator
      number of a $p$-group,  $p$-supersoluble group.
\end{minipage}\end{center}

\vspace{0.5cm}

\section{Introduction }

\indent\\

All groups considered in this paper are finite. We use conventional
notions and notation, as in Huppert \cite{H}.  $G$ always denotes a
finite group,  $|G|$ is the order of $G$,  $\pi(G)$ denotes the set
of all primes dividing  $|G|$, $G_p$ is a Sylow $p$-subgroup of $G$
for some $p\in \pi(G)$. \\

  Suppose that $P$ is a $p$-group for some prime
$p$. Let ${\cal M}(P)$ be the set of all maximal subgroups of $P$.   Let $d_p$ be the smallest generator
      number of a $p$-group $P$, i.e., $p^{d_p} = |P/\Phi(P)|$.   In
\cite{LH}),
 ${\cal M}_{d_p}(P)=\lbrace P_1, ..., P_{d_p}\rbrace$ contains  all elements of ${\cal M}(P)$ such that   $$\bigcap_{i=1}^{d_p}
      P_i=\Phi(P).$$

      We know that ${\cal M}_{d_p}(P)$ is a subset of ${\cal M}(P)$ and
      \noindent so  $$|{\cal M}(P)| >> |{\cal M}_{d_p}(P)|.$$
 For example, suppose that $P$ is  an elementary abelian  $p$-group of order $p^7$. Then $ |{\cal M}(P)| = 137257$, but
 $|{\mathcal M}_{d_p}(P)| = 7$.

\par Two subgroups $H$ and $K$ of $G$  are said to be {\it
permutable} if $HK = KH$.    A subgroup $H$ of $G$ is said to be {\it
s-permutable \cite{K} (or s-quasinormal, $\pi$-quasinormal)} in $G$
if $H$ permutes with every Sylow subgroup of $G$;   Following (\cite{C}), a  subgroup   $H$   of $G$ is {\it semipermutable }  in
$G$  if $H$  permutes with all subgroups $K$  of
$G$   for $ (|H|, |K|) = 1$;  a subgroup $H$  of $G$ is
 {\it s-semipermutable }  in $G$ if $H$ permutes with
every Sylow $p$-subgroup $G_p$ of $G$ with $(|H|, p) = 1$.  Once  the notion of s-semipermutable
 subgroup was introduced, it  has become a hot spot of international group theory research.
 According to MathSciNet, well over 70 papers on this  and related topics
were published in the last decade,   see Refs  \cite{C},  \cite{I}, \cite{LHW}, \cite{LL},  \cite{WLW},  \cite{ZW}, etc.
In \cite{Li2019},   the author  makes a summary  on  the relevant research.

   It is a  interesting topic to determine the structure of a group $G$ by giving some conditions on  ${\cal M}(P)$, where $P$ is a  Sylow $p$-subgroup of $G$ for some prime $p\in \pi(G)$. Srinivasan \cite{S}
proved that if   all maximal  subgroups of any  Sylow subgroups of
$G$ are normal (quasinormal, $S$-quasinormal) in $G$, then $G$ is
supersoluble.  The authors    \cite{LL}  extend  Srinivasan's theorem  as follows:  Suppose that every member in  ${\cal M}_{d_p}
(P)$ is  s-semipermutable in  $G$ for any prime   $p$ in $\pi(G)$.  Then $G$ is  supersoluble (\cite[Theorem 3.4]{LL}).    For the purpose of  perfecting
 the research  of s-semipermutablity,   we apply the idea of localization to extend this  result.  All results in \cite{LL} are unified in the following  theorem. \\

       \noindent{\bf Main Result} ~~   {\it Suppose that $G$ is a group and  $p$ is a fixed prime number in $\pi(G)$ and
    $P$ is  a Sylow
$p$-subgroup of $G$.  Suppose that every member in  ${\cal M}_{d_p}
(P)$ is s-semipermutable in
 $G$.  Then either $P$ is of order $p$ or
$G$ is $p$-supersoluble.}\\

\section{Preliminaries}

\vskip 0.1cm

\noindent{\bf Lemma 2.1}\hspace{0.5cm}(\cite[Lemma 2.1]{LQSW}){\it Let $G$ be a group.

(1) An s-permutable subgroup of $G$ is subnormal in $G$;

(2) If $H \leq K \leq G$ and $H$ is s-permutable in $G$, then $H$ is
s-permutable in $K$;

(3) If $H$ is s-permutable Hall subgroup of  $G$,  then $H
\triangleleft G$;

(4) Let $K \triangleleft G$ and $K \leq H$. Then $H$ is s-permutable
in $G$ if and only if $H/K$ is s-permutable in $G/K$.

(5)  If  $H$, $K$ are s-permutable in $G$, then   $H\cap K$ is also
s-permutable in $G$.

(6) Suppose that $P$ is a  $p$-subgroup of $G$ for some prime $p$.
Then $P$ is s-permutable in $G$ if and only if  $N_G(P) \geq
O^p(G)$.}\\

\par {\bf Lemma 2.2}~~(\cite[Lemma 2.2]{LQSW}) {\it Let $G$ be a group.  Suppose that $H$ is an  s-semipermutable subgroup of
$G$.  Then

(1) If $H \leq K \leq G$, then $H$ is s-semipermutable in $K$.

(2) Let $N$ be a normal subgroup of $G$.  If  $H$ is a $p$-group for
some prime $p\in \pi(G)$, then $HN/N$ is s-semipermutable in $G/N$;

(3) If $ H\leq O_p(G)$, then $H$ is  s-permutable  in $G$.

(4)  Suppose that $H$ is a  $p$-subgroup of $G$ for some prime $p\in
\pi(G)$ and $N$ is  normal in $G$. Then $H\cap N$ is also an
s-semipermutable subgroup of $G$.}\\

The following  useful result is  due to Isaacs (\cite{I}).\\

 \noindent {\bf Lemma 2.3}~~{\it
   Suppose that $H$ is a  $p$-subgroup of $G$. If  $H$ is  s-semipermutable in $G$, then $H^G$ is soluble.}\\

  \noindent{\bf Lemma 2.4} ~~ (\cite[I, Hauptsatz 17.4]{H})  \ {\it Suppose that
$N$ is an  abelian normal subgroup of $G$ and $N\leq M \leq G$ such
that $(|N|, [G:M]) = 1$. If $N$ is complemented in $M$, then $N$ is
complemented in $G$.}\\

\noindent{ \bf Lemma 2.5} ~~(\cite[Lemma 2.4]{ELL})~~{\it  Suppose that $H$ is a nonabelian
simple group. If the Sylow $p$-subgroups $H_p$ of $H$ are  of order
$p$, where  $p$ is  a prime, then the  out automorphism group
$Out(H)$ of $H$ is a  $p^\prime$-group.}

  \section{The proof of main result}

 Suppose that the theorem is false and $G$ is a
 counter-example with minimal order.
  We will derive a  contradiction in several steps.\\

  Step 1.  $O_{p'}(G)=1$.

Denote $N=O_{p'}(G)$. If $N>1$, we consider the factor group $G/N$.
Obviously, $PN/N$ is a Sylow $p$-subgroup of $G/N$, which is
isomorphic   to $P$, so $PN/N$ has the same  smallest generator
number  as $P$, i.e., $d_p$ and so
       $${\cal M}_{d_p}(P/N)=\lbrace P_1/N, ..., P_{d_p}/N\rbrace.$$
        We  know that every  $P_i/N$   is also a
      s-semipermutable subgroup of  $G/N$ by  Lemma 2.2.
                 Thus $G/N$  satisfies the hypotheses of the theorem. We have that either $PN/N$ is of order $p$ or
           $G/O_{p'}(G)$ is $p$-supersoluble
           by the choice of $G$. It follows that
 either $P$ is of order $p$ or
           $G$ is $p$-supersoluble,  a contradiction.
           Thus, we have $N=O_{p'}(G)=1$, as desired.

\vskip 0.15cm

Step 2. $P$ is non-cyclic.

    If $P$ is cyclic, then the unique maximal subgroup $\Phi(P)$ of $P$ is
   semipermutable  in $G$ by the hypotheses. Hence either $P$ is of order $p$ or $G$ is
$p$-supersoluble by  \cite[Theorem 3.2]{LQSW}, a contradiction. \\

Step 3.  $\Phi(P)_{G}=1$. Therefore, $O_{p}(G)$ is an elementary
abelian group.

If not, take any $T\leq \Phi(P)_{G}$ such that $T\unlhd G$. We
        consider the factor group $G/T$.  Since every maximal subgroup of $P$ contains $\Phi(P)$ and
       $P/T$ has the  same smallest generator number  as $P$, so
       $${\cal M}_{d_p}(P/T)=\lbrace P_1/T, ..., P_{d_p}/T\rbrace.$$
        We  know that every  $P_i/T$   is also a
      s-semipermutable subgroup of  $G/N$ by  Lemma 2.2.
            Thus, $G/T$ satisfies the hypotheses of the theorem.
           Hence, either $P/T$ is of order $p$ or $G/T$ is $p$-supersoluble by the
choice of $G$.  If $P/T$ is of order $p$, then $P$ is cyclic,
contrary to Step 2. Hence $G/T$ is $p$-supersoluble. Then $G$ is
$p$-supersoluble, a   contradiction.\\

Step 4.   If $N$ is minimal normal in $G$ contained in $P$, then
 $|N| = p$.

      If $N \leq
  P_i$ for all $P_i\in {\cal M}_{d_p}(P)$,  then
  $$N \leq \cap_{i=i}^{d_p}P_i =
   \Phi(P),$$
   which is contrary to Step 3.  Hence
   there exists a $P_{i_0}\in {\cal
  M}_{d_p}(P)$ such that  $N\not\leq P_{i_0}$.  Then $P_{i_0}\cap N$ is s-semipermutable in $G$ by Lemma 2.2(4). Thus
  $P_{i_0}\cap N$ is s-permutable in $G$ by Lemma 2.2(3). Since $P_{i_0}\cap N$ is  normal in $P$. Hence $P_{i_0}\cap N$ is  normal in $G$ by Lemma 2.1(4). Hence $P_{i_0}\cap N = 1$.    Since $P_{i_0}$ is maximal
  in $P$, we have $N$ is   of order $p$.\\

       Step 5.  All  minimal normal subgroups of $G$ are contained in
   $O_p(G)$.

   Assume that $N$ is a minimal normal subgroup of $G$ which is not a
   $p$-subgroup. As $O_{p^\prime}(G) = 1$ by Step 1, we have  that
   $p||N|$ and $N=N_1\times N_2\times\cdots \times N_s$ is  non-abelian characteristic  simple group, where all $N_i$ are conjugated non-abelian
   simple groups.\\

   (5.1)  $P_i\cap N = 1$ for any $P_i$  in  ${\cal M}_{d_p}(P)$ and  $N$ is a non-abelian simple group, i.e., $s= 1$.

   Pick arbitrary $P_i$ in  ${\cal M}_{d_p}(P)$.  By Lemma 2.2, $P_i\cap N$ is   an  s-semipermutable $p$-subgroup of $G$ by hypotheses and Lemma 2.2(4).  Hence  $P_i\cap N$ is   an  s-semipermutable $p$-subgroup of $N$ by Lemma 2.2(1).  Then $(P_i\cap N)^N$ is  soluble by Lemma 2.3. This implies that   $P_i\cap N = 1$.  Thus    $|P\cap N|\leq p$.  Hence  $N = N_1$ is a simple group. \\

   (5.2)   $O_p(G) =1$.

  If $O_p(G) \not=1$, we can pick a minimal normal subgroup $H$ of
  $G$ contained in $O_p(G)$. By Step 4 we know that $H$ is of order
  $p$.    Hence $HN = H\times N$.

  For any $P_i$  in  ${\cal M}_{d_p}(P)$, if $H\cap P_i = 1$, then $HN \cap P_i = (H\cap P_i)(N\cap P_i) = 1$.   Calculating
  the order of $HNP_i$, this is a contradiction.
  Hence $H\cap P_i \not= 1$. Then $H\leq P_i$, then $H\leq \Phi(P)$, a contradiction. \\

  (5.3)  $C_G(N) =1$.

  Suppose that  $C_G(N) \not= 1$.  Now  we  pick a minimal normal
    subgroup $N^*$ of $G$ contained in $C_G(N)$. Then $N^*$  is  non-abelian simple  by  (5.1) and (5.2) and  $N\cap N^* = 1$.

  Since  $P_1\cap NN^* = (P_1\cap N)(P_1\cap N^*) = 1$, we have  $|P_1NN^*|_p =    p^2|P_1| >|P|$, a contradiction. \\

      (5.4) Finishing the proof of Step 5.

   By (5.3) we have    $C_G(N)=1$. Then $G$ and $G/N$ are isomorphic to a
subgroup of $Aut(N)$ and a subgroup of $Aut(N)/Inn(N)$,
respectively. This means that $N_p$ is of order $p$ and $p$ divides
the order  of $Out(N)$.
By Lemma 2.5, this is impossible. \\

   Step 6. $G = O_p(G)\rtimes M$, the semi-direct product of
  $O_p(G)$ with a subgroup $M$ of $G$ and $O_p(G)$ is a direct
  product of normal subgroups of $G$ of order $p$.

  Let $N_1$ be a minimal normal subgroup of $G$ contained in $O_p(G)$.
   Then $N_1$ is of order $p$ by Step 4  and $N_1\cap \Phi(P) =1$ by Step 3.
   Hence there exists a maximal subgroup $S_1$ of $P$
    such that $N_1 \cap S_1 = 1$. By Lemma 2.4,
   $N_1$ has a complement  $K$ in $G$, i.e., $G = N_1K$ and $N_1\cap
   K = 1$. Then $O_p(G) = N_1(O_p(G)\cap K)$. It is easy to see that
   $O_p(G)\cap K$ is normal in $G$ and $P\cap K$ is a Sylow $p$-subgroup of $K$.
   If $O_p(G)\cap K = 1$, then Step 6  holds. So assume that $O_p(G)\cap K\not= 1$. Then we can
  pick a minimal subgroup $N_2$ contained in $O_p(G)\cap K$. By Step 4, $N_2$ is of order $p$ and there exists a  maximal subgroup $S_2$  of $P$   such that $N_2 \cap S_2 = 1$. Then $P = N_2S_2 = S_2(O_p(G)\cap K) = S_2(P\cap
   K)$. Since $|(P\cap K):(S_2\cap K)| = |S_2(P\cap K):S_2| = |P:S_2| =
   p$,  $S_2\cap K$ is a complement of $N_2$ in $P\cap K$. Therefore
   $N_2$ has a complement $L$ in $K$ by Lemma 2.4. Then $G = N_1K =
   (N_1\times N_2)\rtimes L$. Continuing this process, we have
   finally  $G = O_p(G)\rtimes M$ and $O_p(G) = N_1\times
   N_2\times \cdots\times N_r$, where $N_i$ is a normal subgroup of $G$ of order
   $p$. \\

   Step 7. The final contradiction.

      Since $N\leq Z(P)$ for any minimal normal subgroup $N$ of $G$,
     $P\leq C_G(O_p(G))$.    Since $C_G(O_p(G))\cap M \lhd \langle
     O_p(G), M\rangle = G$, $C_G(O_p(G))\cap M = 1$ by Step 4 and 5. Then
     $P\cap M = 1$. This implies that $P = P\cap O_p(G)M =
     O_p(G)(P\cap M) = O_p(G)$. Therefore by  Step 6 we have that $G$ is $p$-supersoluble, the final
    contradiction. \hfill $\Box$\\

 \noindent{\bf Remark}~~ The author does not know   the proof without using   the  classification  theorem of finite  non-abelian simple groups.\\

\section{Applications}

We give some applications of our Main Result.

 Suppose that
$p$ is
  the smallest prime dividing the order of $G$. We
 know that  $G$ is $p$-nilpotent
 if $G_p$ is cyclic by \cite[IV Satz 2.8]{H} and $p$-supersolubility
 implies the $p$-nilpotency.  By  our Main Result, immediately
  the following corollaries hold.\\

  \noindent {\bf Corollary  4.1}~~(\cite[Theorem 3.1]{LL})~~ {\it Let $p$ be the smallest prime dividing
     $|G|$ and $P$ a Sylow $p$-subgroup of $G$.  Then $G$ is $p$-nilpotent  if and only if
      every member in  ${\cal M}_{d_p}
(P)$ is s-semipermutable in  $G$.}\\

If $G_p$ is   of order $p$ and $G$  is $p$-soluble, then,  obviously,  $G$ is  $p$-supersoluble.  Hence we have:\\

\noindent{\bf  Corollary 4.2} ~~ (\cite[Theorem 3.8]{LL}) {\it Suppose that $G$ is a
$p$-soluble group, where  $p$ is a fixed prime number in $\pi(G)$,
and $P$ is a Sylow $p$-subgroup of $G$. Then $G$ is
$p$-supersoluble if  every member in  ${\cal M}_{d_p}
(P)$ is s-semipermutable  in  $G$.}\\

\noindent{\bf  Corollary 4.3}~~ (\cite[Theorem 3.9]{LL}) {\it Suppose that  $P$ is a Sylow
$p$-subgroup of $G$ and $N_G(P)$ is $p$-nilpotent for some prime
$p\in \pi(G)$. Then $G$ is $p$-nilpotent if  every
member in  ${\cal M}_{d_p}
(P)$ is s-semipermutable in $G$.}\\

    \noindent  {\bf Proof.}~~ Applying  our Main Result we know that either $P$
      is cyclic  or $G$ is
  $p$-supersoluble.
  If $P$ is cyclic, then we have $N_G(P) = C_G(P)$. Applying
  Burnside's $p$-nilpotence criterion (\cite[Hauptsatz IV.2.6]{H}), we get that $G$ is $p$-nilpotent.
  Now suppose that $G$ is $p$-supersoluble.
  Since the $p$-length of $p$-supersoluble groups is at most 1, we
  have $PO_{p^\prime}(G) $ is normal in $G$. Set $\overline{G} = G/O_{p^\prime}(G)$.
  Then  $\overline{G} =
  N_{\overline{G}}(\overline{P}) = N_G(P)O_{p^\prime}(G)/O_{p^\prime}(G)$
  is $p$-nilpotent by hypothesis. Hence $G$ is $p$-nilpotent, as desired. \hfill$\Box$\\

\end{document}